\documentclass[a4paper,leqno]{article}

\usepackage{RRA4}
\usepackage{hyperref}

\RRdate{Juin 2006}
\RRNo{5940}

\RRauthor{
  Hend {Ben Ameur}%
  \thanks{LAMSIN, ENIT, BP 37, Le Belvédère, 1002 Tunis, Tunisie.
    {\tt hbenameur@yahoo.ca}.}
  \and François Clément%
  \thanks{Projet Estime. {\tt Francois.Clement@inria.fr}.}
  \and Pierre Weis%
  \thanks{Projet Cristal. {\tt Pierre.Weis@inria.fr}.}
  \and Guy Chavent%
  \thanks{Projet Estime. {\tt Guy.Chavent@inria.fr}.}
  }
\authorhead{Ben Ameur, Clément, Weis \& Chavent}

\RRtitle{L'algorithme des indicateurs de raffinement multidimensionel pour
  une paramétrisation optimale}
\RRetitle{The Multidimensional Refinement Indicators Algorithm for Optimal
  Parameterization}
\titlehead{Multidimensional Refinement Indicators Algorithm} 

\RRresume{% $Id: resume.tex,v 1.6 2006/06/27 09:01:13 fclement Exp $

L'estimation de paramètres distribués dans des équations aux dérivées
partielles (EDP) à partir de mesures de la solution de l'EDP peut mener à des
problèmes de sous-détermination.
Le choix d'une paramétrisation est un moyen usuel pour ajouter de
l'information a priori en réduisant le nombre d'inconnues en relation avec la
physique du problème.
L'algorithme des indicateurs de raffinement fourni une technique de
paramétrisation adaptative fructueuse qui ouvre parcimonieusement les degrés de
liberté de façon itérative.
Nous présentons une
{\em nouvelle forme générale de l'algorithme des indicateurs de raffinement}
qui s'applique à l'estimation des paramètres multi-dimensionnels dans toute EDP.
Dans le cas linéaire, nous établissons le lien entre l'indicateur de
raffinement et la décroissance de la fonction objectif des moindres carrés
quantifiant l'erreur aux données.
Nous donnons des résultats numériques pour le cas simple du modèle identité, et
cette application permet de voir l'{\em algorithme des indicateurs de
  raffinement comme une technique de segmentation d'image}.
}
\RRabstract{% $Id: abstract.tex,v 1.2 2008/01/16 12:32:01 fclement Exp $

The estimation of distributed parameters in a partial differential equation
(PDE) from measures of the solution of the PDE may lead to underdetermination
problems.
The choice of a parameterization is a frequently used way of adding a priori
information by reducing the number of unknowns according to the physics of the
problem.
The refinement indicators algorithm provides a fruitful adaptive
parameterization technique that parsimoniously opens the degrees of freedom in
an iterative way.
We present a {\em new general form of the refinement indicators algorithm} that
is applicable to the
{\em estimation of distributed multidimensional parameters} in any PDE.
In the linear case, we state the relationship between the refinement indicator
and the decrease of the usual least-squares data misfit objective function.
We give numerical results in the simple case of the identity model, and this
application reveals the {\em refinement indicators algorithm as an image
segmentation technique}.
}

\RRmotcle{problème inverse, paramétrisation, raffinement optimal, segmentation
  d'image, programmation fonctionnelle}
\RRkeyword{inverse problem, parameterization, optimal refinement, image
  segmentation, fonctional programming}

\RRprojets{Estime et Cristal}

\RRtheme{\THNum \THSym}

\URRocq

\usepackage[T1]{fontenc}
\usepackage[latin1]{inputenc}
\usepackage{pslatex}
\usepackage{amssymb}
\usepackage{amsmath}
\usepackage{graphicx}
\usepackage{latexsym}
\usepackage{epic,eepic}
\usepackage{alltt}

\renewcommand{\arraystretch}{2}

\newtheorem{remark}{Remark}

\newcommand{\Id}{{\rm Id}}
\newcommand{\ad}{{\rm ad}}
\newcommand{\init}{{\rm init}}
\newcommand{\meas}{{\rm meas}}
\newcommand{\opt}{{\rm opt}}
\newcommand{\true}{{\rm true}}
\newcommand{\sgn}{{\rm sgn}}
\newcommand{\best}{{\rm best}}
\newcommand{\dichotomy}{{\rm dichotomy}}
\newcommand{\card}{{\rm card}}
\newcommand{\cut}{{\rm cut}}

\newcommand{\R}{{\mathbb R}}

\newcommand{\calF}{{\mathcal F}}
\newcommand{\calI}{{\mathcal I}}
\newcommand{\calJ}{{\mathcal J}}
\newcommand{\calL}{{\mathcal L}}
\newcommand{\calM}{{\mathcal M}}
\newcommand{\calO}{{\mathcal O}}
\newcommand{\calP}{{\mathcal P}}
\newcommand{\calT}{{\mathcal T}}

\newcommand{\programminglanguage}[1]{{\sc #1}}

\newcommand{\ocaml}{\programminglanguage{OCaml}}
\newcommand{\ocamlpppl}{\programminglanguage{OCamlP3l}}
\newcommand{\clang}{\programminglanguage{C}}
\newcommand{\cpp}{\programminglanguage{C$++$}}
\newcommand{\fortran}{\programminglanguage{Fortran}}

\begin{document}
\makeRR

% $Id: introduction.tex,v 1.2 2008/01/16 12:32:01 fclement Exp $

\section{Introduction}

The inverse problems of parameter estimation consist in identifying unknown
parameters which are space dependent coefficients in a system of partial
differential equations modeling a physical problem.
The parameter estimation problem can be set as a minimization problem of an
objective function defined as the least-squares misfit between measurements and
the corresponding quantities computed with a chosen parameterization of the
parameters (solution of the partial differential equations system under
consideration).

Experimental measurements are expensive, thus one of the difficulties when
solving a parameter estimation inverse problem is that the data is usually
insufficient to estimate the value of the parameter in each cell of the
computational mesh, the parameter estimation problem is underdetermined.
Therefore, a parameterization of the sought parameter is chosen to reduce the
number of unknowns.
The multiscale parameterization approach is popular for many applications,
e.g. see~\cite{liu1,bunsalzalcha,prattetal,akcbirgha}.
It consists in solving the problem through successive approximations by
refining uniformly the parameter at the next finer scale all over the domain
and stopping the process when the refinement does not induce significant
decrease of the objective function any more.
This method is very attractive as it may provide a regularizing effect on the
resolution of the inverse problem and it may also circumvent local minima
problems as shown in~\cite{ch5,clechagom}.
But its main drawback is the over-parameterization problem due to the uniform
refinement of the parameter, see~\cite{chacha,mose}.

To avoid this pitfall, the refinement indicators algorithm provides an
adaptative parameterization technique that parsimoniously opens the degrees of
freedom in an iterative way driven at first order by the model to locate the
discontinuities of the sought parameter.
The detailed definition of the technique was described in~\cite{benchajaf} for
the application to the estimation of hydraulic transmissivities, and a variant
of the idea was first briefly presented in~\cite{chabis} for the application to
the history matching of an oil reservoir.

In the present work, we extend the definition of refinement indicators to the
more general case of multidimensional distributed parameters, we make the link
between the (first order) refinement indicator and the exact decrease of the
objective function in the linear case, and we propose a generic version of the
algorithm that may be applied to any inverse problem of parameter estimation
in a system of partial differential equations.
We show numerical results for the simple case of the identity model.
The application of the refinement indicators algorithm to the identity model
applied to images is interpreted as an amazing image segmentation technique.

The paper is organized as follows.
Section~\ref{sec:indicators} is devoted to theoretical developments about
adaptive parameterization through the refinement indicators notion: we quantify
the effect of refinement on the decrease of the objective function, and we
generalize the definition of refinement indicators for multidimensional
parameters.
We present in Section~\ref{sec:algorithm} the general form of the refinement
indicators algorithm, we discuss the issue of finding the best cutting for
multidimensional parameters and briefly point out implementation subtleties
of the algorithm.
In section~\ref{sec:segmentation}, we focus on the identity model; we then
apply the multidimensional refinement indicators algorithm to the particular
case of RGB color images.

% $Id: indicators.tex,v 1.2 2008/01/16 12:32:01 fclement Exp $

\section{Adaptive inversion and refinement indicators}
\label{sec:indicators}

The inverse problem of parameter estimation is first set as the minimization of
a least-squares objective function.
Then, the notion of adaptive parameterization induced by zonations of the
domain is presented in the discrete case for piecewise constant parameters.
The decrease of the optimal value of the objective function resulting from the
splitting of one zone into two subzones is then quantified for the linearized
problem in the continuous case, and refinement indicators are finally defined
in the general case from a first order approximation.

\subsection{Least-squares inversion}
\label{sec:ls}

When defining a system of partial differential equations (PDE) modeling a
physical behavior in a domain $\Omega\subset\R^N$ (e.g. with $N=1$, 2 or 3),
one usually introduces parameters such as material properties.
We consider the case where these parameters are \emph{distributed} and possibly
\emph{vector valued}, i.e. belong to a space~$P$ of functions defined
over~$\Omega$ with values in $\R^{n_p}$ ($n_p\geq 1$ is the \emph{dimension} of
the vector parameter $p(x)$ at any $x\in\Omega$).

Some of these parameters may not be accessible by direct measurements, so they
have to be estimated indirectly.
If $p_\true\in P$ denotes such an unknown parameter that we are looking for in
a set of admissible parameters~$P^\ad$, and if
$d\simeq\calF(p_\true)\in\R^{n_d}$ (the \emph{data}) denotes the corresponding
vector of available measurements of the solution of the PDE, one can attempt an
indirect determination of~$p_\true$ from~$d$ by solving the
\emph{least-squares inverse problem}:
\begin{equation}
  \label{eqn:lsp}
  \mbox{minimize } J(p) \mbox{ for } p\in P^\ad
\end{equation}
where $J(p)$ is the least-squares misfit between~$d$ and the corresponding
quantities $\calF(p)$ computed from the current parameter~$p$,
\begin{equation}
  \label{eqn:J}
  \calJ(p) = \frac{1}{2} \| d - \calF(p) \|^2.
\end{equation}

The \emph{forward operator}~$\calF$ is the composition of the
\emph{model operator}~$\calM$ (which computes the solution of the PDE for a
given parameter~$p$), with the \emph{observation operator}~$\calO$ (which
computes the output of the observation device applied to the solution of the
PDE).
For example, observations can be made of the values of the solution of the PDE
at a set of measurement points of the domain~$\Omega$.

Minimizing such an objective function avoids the difficulty of
inverting~$\calM$, which is impossible in most cases.

\subsection{Adaptive parameterization}
\label{sec:adapt}

In parameter estimation problems, the data is usually insufficient to estimate
the value of the (possibly vector valued) unknown parameter at each point
$x\in\Omega$ (i.e. in each grid cell after discretization).
It is usually impossible to increase the number of measurements, and one
solution is then to reduce the number of unknowns by searching for the
parameter in a subspace of~$P$ of finite (and small) dimension.

For that, we proceed in two steps: first we construct a finite dimensional
subspace~$P_n$ of the infinite dimensional space~$P$, then we look for an
approximation~$p_n$ of the unknown and infinite dimensional true parameter in
$P_n^\ad=P_n\cap P^\ad$.
The dimension of~$P_n$ should remain small enough to avoid
over-parameterization, i.e. the Jacobian of~$\calF$ as a function from~$P_n$
to~$\R^{n_d}$ should have a full rank, and the approximation~$p_n$ should
explain the data up to the noise level.
Notice that here both~$P_n$ and~$p_n$ are unknown.

For a given subspace~$P_n$, a natural choice for~$p_n$ is a minimizer of the
objective function~$J$ in $P_n^\ad$.

It is classical to build the subspace~$P_n$ in an iterative way by successive
refinements as in the multiscale approach,
see~\cite{liu1,bunsalzalcha,prattetal,akcbirgha}.
We use the \emph{adaptive parameterization} technique developed
in~\cite{benchajaf}.
These approaches are also known to have a regularizing effect on the resolution
of the inverse problem, e.g. see~\cite{ch5,liu1,liu2,engl,benkal}.

The adaptive parameterization approach consists in adding one (vectorial)
degree of freedom at a time.
The new degree of freedom is chosen according to \emph{refinement indicators}
which evaluate the benefit for the objective function of the introduction of a
given family of degrees of freedom.
Hence, each subspace~$P_n$ is of dimension~$n\times n_p$.
The process is stopped when $p_n$ explains the data up to the noise level.

It is convenient to consider~$P_n$ as the range of
an---unknown---\emph{parameterization map}
\begin{equation}
  \label{eqn:P_n}
  \calP_n: m_n \in M_n^\ad \longmapsto p_n \in P_n^\ad
\end{equation}
where~$m_n$ is the \emph{coarse} parameter (of small finite dimension), by
opposition to the \emph{fine} parameter $p_n$ (of large, and possibly infinite
dimension).
$M_n^\ad$ is the space of admissible coarse parameters.
Typically, $m_n$ is made of the coefficients of $p_n\in P_n$ on a basis of
$P_n$, in which case $\calP_n$ is a linear operator.
But the parameterization map can also be nonlinear, as for example in the case
where~$m_n$ is made of the coefficients of a closed form formula defining the
fine parameter $p_n$.
For any parameterization map~$\calP_n$, we define the same objective function
on~$M_n^\ad$ by
\begin{equation}
  \label{eqn:J_n}
  \calJ_n(m_n)=\calJ(\calP_n(m_n)),
\end{equation}
and the least-squares problem becomes:
\begin{equation}
  \label{eqn:lspn}
  \mbox{minimize } J_n(m_n) \mbox{ for } m_n\in M_n^\ad.
\end{equation}

Let $Z_n=(Z_{n,j})_{1\leq j\leq n}$ be a partition of the closure of the
domain~$\Omega$ in~$n$ parts, i.e. for all $j\in\{1,\ldots,n\}$
$Z_{n,j}\subset\overline{\Omega}$,
$\bigcup_{1\leq j\leq n}Z_{n,j}=\overline{\Omega}$, and
$Z_{n,j}\cap Z_{n,k}=\emptyset$ as soon as $j\not=k$.
We consider that the subspace~$P_n$ is made of piecewise constant
(vector valued) functions on this partition.
We call \emph{zonation} the partition~$Z_n$, and \emph{zone} each part~$Z_{n,j}$
of the zonation.
The parameterization map associated with the zonation~$Z_n$ is then
\begin{equation}
  \label{eqn:P_n2}
  \calP_n: m_n = (m_{n,j})_{1\leq j\leq n}
  \quad \longmapsto \quad
  p_n =
  (x \in \Omega \mapsto  m_{n,j} \in \R^{n_p} \mbox{ when } x \in Z_{n,j}).
\end{equation}
It associates the coarse parameter $m_n$ with the function $p_n$ which takes on
each zone $Z_{n,j}$ the constant value $m_{n,j}$.

In the iterative process, when going from~$P_n$ to~$P_{n+1}$, we introduce a
new (vectorial) degree of freedom by dividing one zone of~$Z_n$ into two
subzones, thus producing a new zonation~$Z_{n+1}$ having one zone more
than~$Z_n$.

\subsection{Quantifying the effect of refinement for the linearized problem}
\label{sec:quantif}

% cf. Trond Mannseth (for mesh refinement) use a more accurate indicator, but
% also more expensive (solve the linearized problem)

Considering a current zonation, the best refinement would be the one
corresponding to the largest decrease of the objective function.
This decrease can only be computed by actually solving the least-squares
problem~(\ref{eqn:lspn}) for each possible refinement to be tested, hence only a
very small number of them could be evaluated before reaching a decision.
So, we search for a closed form formula for the decrease of the objective
function for the \emph{linearized} problem.
We shall still denote by $\calF$ the linearized forward map.

For simplicity, we consider the case where the optimization
problem~(\ref{eqn:lspn}) is  unconstrained, so that $P_n^\ad=P_n$ and
$M_n^\ad=\R^{n \times n_p}$ for all~$n$.

Let us compute the decrease of the optimal objective function~(\ref{eqn:J_n})
resulting from refining the current zonation by splitting one of the zones into
two subzones by means of a continuum of intermediate minimization problems
constraining the discontinuity jump at the interface between the subzones.

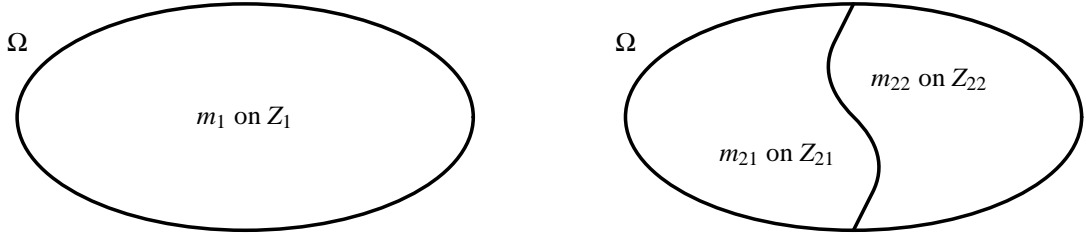
\begin{figure}[ht]
  \label{fig:P1P2}
  \begin{center}
    % $Id: cutting.tex,v 1.2 2006/05/23 16:43:23 fclement Exp $

\setlength{\unitlength}{1mm}

\begin{picture}(140,30)(0,0)
  \allinethickness{0.5mm}

  \put(30,15){\ellipse{60}{30}}
  \put(0,25){\makebox(0,0){$\Omega$}}
  \put(30,15){\makebox(0,0){$m_1$ on $Z_1$}}

  \put(110,15){\ellipse{60}{30}}
  \put(110,0){\spline(0,0)(5,10)(-5,20)(0,30)}
  \put(80,25){\makebox(0,0){$\Omega$}}
  \put(100,10){\makebox(0,0){$m_{21}$ on $Z_{21}$}}
  \put(120,20){\makebox(0,0){$m_{22}$ on $Z_{22}$}}
\end{picture}
  \end{center}
  \caption{Refinement of a single zone zonation represented in a
    two-dimensional space for simplicity.
    Left: one-zone parameterization~$\calP_1$ with $\overline{\Omega}=Z_{1,1}$.
    Right: two-zone parameterization~$\calP_2$ with
    $\overline{\Omega}=Z_{2,1}\cup Z_{2,2}$.}
\end{figure}

Let us consider the case where only one zone covers the whole domain,
corresponding to the one-zone parameterization~$\calP_1$, see
Figure~\ref{fig:P1P2} (left).
The refinement builds the two-zone parameterization~$\calP_2$
of~Figure~\ref{fig:P1P2} (right).
We denote by~$m_1^\opt=(m_{1,1}^\opt)$ and~$\calJ_1^\opt$,
$m_2^\opt=(m_{2,1}^\opt,m_{2,2}^\opt)^T$ and~$\calJ_2^\opt$ the optimal coarse
parameter and objective function respectively when considering the
parameterizations~$\calP_1$ and~$\calP_2$.

If the discontinuity jump $c^\opt=m_{2,1}^\opt-m_{2,2}^\opt$ was known, then
minimizing~$\calJ_2$ (i.e. minimizing~$\calJ$ considering the
parameterization~$\calP_2$) under the constraint $m_{2,1}-m_{2,2}=c^\opt$
would give us the same optimal coarse parameter~$m_2^\opt$ and the same optimal
objective function value~$\calJ_2^\opt$, the ones obtained without any
constraints.
And minimizing~$\calJ_2$ and under the constraint $m_{2,1}-m_{2,2}=0$
would keep the optimal values~$m_1^\opt$ and~$\calJ_1^\opt$, obtained with the
parameterization~$\calP_1$.
Thus, when the discontinuity jump~$c$ goes continuously from~$0$ to~$c^\opt$,
then the minimum of~$\calJ_2$ under the constraint $m_{2,1}-m_{2,2}=c$
goes continuously (and even in a continuously differentiable way when
$\calF$~is continuously differentiable) from the minimum obtained for a single
zone to the minimum obtained with the two zones.

The components~$m_{2,1}$ and~$m_{2,2}$ of the coarse parameter are column
vectors of dimension $n_p$ (i.e. $n_p=1$ in the scalar case).
We can denote the constraint on the discontinuity in matrix form by
\begin{equation}
  \label{eqn:constraint}
  A m_2 = c
\end{equation}
where $m_2=(m_{2,1},m_{2,2})^T$ and the $n_p\times 2n_p$ rectangular matrix~$A$
is of the form
\begin{equation}
  \label{eqn:A}
  \renewcommand{\arraystretch}{1}
  A =
  \left(
  \begin{array}{rrrrrrrr}
    1      & 0      & \cdots & 0      & -1     & 0      & \cdots & 0      \\
    0      & \ddots & \ddots & \vdots & 0      & \ddots & \ddots & \vdots \\
    \vdots & \ddots & \ddots & 0      & \vdots & \ddots & \ddots & 0      \\
    0      & \cdots & 0      & 1      & 0      & \cdots & 0      & -1
  \end{array}
  \right).
\end{equation}

We denote the (linear) direct operator in terms of the coarse parameter~$m_n$ by
$\calF_n=\calF\circ\calP_n$.
In the following computation of the decrease of the optimal value of the
objective function, we will only consider the two-zone
parameterization~$\calP_2$.
Hence, when no ambiguity can arise, we will simply write~$\calF$ for~$\calF_2$,
$\calJ$ for~$\calJ_2$ and~$m$ for~$m_2$.
The gradient of the quadratic objective function $\calJ(=\calJ_2)$ is given by
\begin{equation}
  \label{eqn:GJP}
  \nabla\calJ(m) = \calF^T (\calF m - d).
\end{equation}

The Lagrangian function associated with the minimization of~$\calJ$ under the
constraint~(\ref{eqn:constraint}) writes
\begin{equation}
  \label{eqn:Lc}
  \calL^c(m,\lambda) = \calJ(m) + \left< \lambda, A m - c \right>
\end{equation}
where~$\lambda$ is the Lagrange multiplier associated with the constraint.
It is a vector of dimension~$n_p$.
Then, the Lagrange condition ensures that at the optimum,
$m^{c,\opt}=\arg\min_{Am=c}\calJ(m)$ and~$\lambda^{c,\opt}$ are obtained by
solving
\begin{equation}
  \label{eqn:lagrcond}
  \left\{
  \begin{array}{l}
    \displaystyle
    \frac{\partial\calL^c}{\partial m}(m,\lambda) =
    \nabla\calJ(m) + A^T \lambda = 0, \\
    \displaystyle
    \frac{\partial\calL^c}{\partial\lambda}(m,\lambda) = A m - c = 0.
  \end{array}
  \right.
\end{equation}
Using the expression~(\ref{eqn:GJP}) of the gradient, the Lagrange condition
for $c=0$ becomes
\begin{equation}
  \label{eqn:lagrcond0}
  \left\{
  \begin{array}{l}
    \calF^T (\calF m^{0,\opt} - d) + A^T \lambda^{0,\opt} = 0, \\
    A m^{0,\opt} = 0.
  \end{array}
  \right.
\end{equation}
The last equation writes $m_{2,1}^{0,\opt}=m_{2,2}^{0,\opt}(=m_{1,1}^\opt)$,
and we have
\begin{equation}
  \label{eqn:J20opt}
  \calJ^{0,\opt} = \calJ_{\calP_2}(m^{0,\opt}) = \calJ_1^\opt.
\end{equation}

The optimal coarse parameter $m^\opt(=m_2^\opt)=\arg\min\calJ(m)$ (i.e., without
constraint) satisfies the optimality condition
\begin{equation}
  \label{eqn:optcond}
  \nabla\calJ(m^\opt) = \calF^T (\calF m^\opt - d) = 0.
\end{equation}

Let $e=m^\opt-m^{0,\opt}$.
Developing the squared norm yields
\begin{equation*}
  \calJ^\opt = \calJ^{0,\opt}
  + \frac{1}{2} \| \calF e \|^2
  + \left< \calF m^{0,\opt} - d, \calF e \right>.
\end{equation*}
But, from the optimality condition~(\ref{eqn:optcond}), we have
\begin{equation*}
  \left< \calF m^{0,\opt} - d, \calF e \right> = - \| \calF e \|^2
\end{equation*}
and taking the difference between~(\ref{eqn:optcond}) and the first equation
of~(\ref{eqn:lagrcond0}) leads to
\begin{equation*}
  \renewcommand{\arraystretch}{1}
  \calF^T \calF e = A^T \lambda^{0,\opt} =
  \left(
  \begin{array}{r}
    \lambda^{0,\opt} \\
    -\lambda^{0,\opt}
  \end{array}
  \right).
\end{equation*}
Therefore, when $\calF^T\calF$ is invertible%
\footnote{E.g., when~$\calF$ is injective in the finite-dimensional case.}, we
finally obtain
\begin{equation}
  \label{eqn:decrease}
  \calJ_1^\opt - \calJ_2^\opt =
  \frac{1}{2} \| A^T \lambda^{0,\opt} \|_{(\calF^T \calF)^{-1}}^2
  = \frac{1}{2} \left< A^T \lambda^{0,\opt},
  (\calF^T \calF)^{-1} A^T \lambda^{0,\opt} \right>
\end{equation}
where the norm and the scalar product are now defined in the space of coarse
parameters for the refined zonation (and not in the space of data as before).

\begin{remark}
  \label{rem:nzone}
  In the general case, an $n$-zone parameterization is refined into an
  $(n+1)$-zone parameterization by splitting only one zone into two subzones.
  Then, the matrix defining the new discontinuity jump of the parameter is a
  block matrix containing only one nonzero block corresponding to the
  discontinuity of the parameter between the two subzones and this block is
  equal to the matrix~$A$ of~(\ref{eqn:A}).
  Hence, the computation of the decrease of the optimal value of the objective
  function is exactly the same.
\end{remark}

Because the number of zones is supposed to stay small, after discretization,
the direct operator~$\calF_n$ is usually injective.
But applying~$(\calF_n^T \calF_n)^{-1}$ is at least as expensive as solving the
linearized minimization problem.
Thus, the exact computation of the decrease of the optimal value of the
objective function for the linearized problem is generally not compatible with
the idea of a very fast computation of an indicator on the quality of the new
zonation.
Nevertheless, an important exception is the case of the identity direct
operator; it will be developed in section~\ref{sec:segmentation}.
The next subsection is devoted to a first order computation in the general
nonlinear case.

\subsection{Refinement indicators for multidimensional parameters}
\label{sec:refind}

Following the notations of the previous subsection, let~$\calJ^{c,\opt}$
denote the optimal value of the objective function obtained with the
parameterization~$\calP_2$ under the constraint~(\ref{eqn:constraint}).
A first order development of~$\calJ^{c,\opt}$ with respect to the discontinuity
jump~$c$ at $c=0$ can be written
\begin{equation*}
  \calJ^{c,\opt} = \calJ^{0,\opt}
  + {\frac{\partial\calJ^{c,\opt}}{\partial c}}_{|c=0}c + \ldots
\end{equation*}

The norm of the quantity
$\displaystyle{\frac{\partial\calJ^{c,\opt}}{\partial c}}_{|c=0}$ tells us, at
the first order, how large would be the difference between~$\calJ^{c,\opt}$
and~$\calJ^{0,\opt}=\calJ_1^\opt$ (from~(\ref{eqn:J20opt})).
So, it gives us the first order effect on the optimal value of the objective
function produced by refining and allowing a discontinuity of the parameter
between the two subzones of the parameterization~$\calP_2$.
Similarly to the particular scalar case developed in~\cite{benchajaf}, this
norm is called the {\em refinement indicator} associated with the splitting of
the zone~$Z_{1,1}$ into the two subzones~$Z_{2,1}$ and~$Z_{2,2}$, and it is
denoted by~$\calI$.

Deriving the expression~(\ref{eqn:Lc}) of the Lagrangian with respect to~$c$
gives
\begin{equation*}
  \frac{\partial\calL^c}{\partial c}(m,\lambda) =
  \frac{\partial\calJ}{\partial c}(m) - \lambda.
\end{equation*}
Hence, from the Lagrange condition~(\ref{eqn:lagrcond}), we have at the
optimum for $c=0$
\begin{equation*}
  \frac{\partial\calJ}{\partial c}(m^{0,\opt}) - \lambda^{0,\opt} = 0.
\end{equation*}
Then, since
\begin{equation*}
  {\frac{\partial\calJ^{c,\opt}}{\partial c}}_{|c=0}
  = {\frac{\partial\calJ}{\partial c}(m^{c,\opt})}_{|c=0}
  = \frac{\partial\calJ}{\partial c}(m^{0,\opt}),
\end{equation*}
this shows that the refinement indicator is nothing but the norm of the
Lagrange multiplier,
\begin{equation}
  \label{eqn:refind}
  \calI = \| \lambda^{0,\opt} \|.
\end{equation}

\begin{remark}
  \label{rem:highorder}
  In this section, we have considered that the vector parameter was a piecewise
  constant function (i.e. constant in each zone).
  The idea of refinement indicators is based on defining a zonation by
  identifying the interfaces between the zones that are related to the
  discontinuities of the parameter.
  To estimate more regular parameters, one can use similar refinement
  indicators in the more general case where the vector parameter is a piecewise
  higher order function, which is continuous in each zone and presents
  discontinuities at the interfaces between the zones.
\end{remark}

% $Id: algorithm.tex,v 1.2 2008/01/16 12:32:01 fclement Exp $

\section{Refinement indicators algorithm}
\label{sec:algorithm}

To solve numerically the parameter estimation problem, we consider a
discretization of the domain~$\Omega$ by a fine mesh $\calT_h=(K_i)_{i\in I}$,
i.e. with $\overline{\Omega}=\bigcup_{i\in I}K_i$.
Let $n_I=\card I$ be the number of cells.
The fine parameter is then the vector $p=(p_i)_{i\in I}$ which approximates the
parameter on each cell~$K_i$ of~$\calT_h$.
We consider zonations~$Z_n=(Z_{n,j})_{1\leq j\leq n}$ following the mesh such
that the associated parameterization map~(\ref{eqn:P_n}) is
\begin{equation}
  \label{eqn:P_h}
  \begin{array}{rrcl}
    \calP_n: & (\R^{n_p})^n & \longrightarrow & (\R^{n_p})^{n_I} \\
    & m_n = (m_{n,j})_{1\leq j\leq n} & \longmapsto & p_n = (p_{n,i})_{i\in I}
    \mbox{ with }  p_{n,i} = m_{n,j} \mbox{ when } K_i \subset Z_{n,j}.
  \end{array}
\end{equation}
In other words, when dealing with a discrete problem, a zonation is related to
a partition of the set~$I$ of indices:
\begin{equation*}
  Z_{n,j} = \bigcup_{i\in I_{n,j}} K_i
  \qquad \mbox{with} \quad
  I = \bigsqcup_{1\leq j\leq n} I_{n,j}
  \mbox{ (disjoint union).}
\end{equation*}

The finite set of measurement points is indexed by $I_m\subset I$.
Keeping the same notations, the discrete objective function to minimize writes
\begin{equation}
  \label{eqn:J_h}
  \calJ(p) = \frac{1}{2} \sum_{i\in I_m} (d_i - u_i)^2
\end{equation}
where $d_i$~is the measurement in the cell~$K_i$ and~$u_i$ denotes
 the discrete approximation of $\calF(p)$ in the same cell~$K_i$.

\subsection{The adaptive parameterization technique}
\label{sec:technique}

In practice, the optimal coarse parameter~$m_n^\opt$ associated with the
current parameterization~$\calP_n$ is computed by applying a gradient algorithm
to the minimization of the corresponding objective
function $\calJ_n$.
Hence, not only $m_n^\opt=(m_{n,j}^\opt)_{1\leq j\leq n}$ is available, but
also the coarse gradient $\nabla\calJ_n(m_n^\opt)$, which vanishes at the
reached minimum.
The key is to compute the gradient of~$\calJ_n$ for the current
parameterization by the adjoint approach considering the fine
discretization~$\calT_h$ of~$\Omega$.
This approach provides, at no additional cost, the \emph{fine gradient}
\begin{equation}
  \label{eqn:GpJn}
  \nabla_p \calJ (p_n^\opt) =
  \left( \frac{\partial\calJ}{\partial p_i}(p_n^\opt) \right)_{i\in I}
\end{equation}
where $p_n^\opt=\calP_nm_n^\opt$, since the \emph{coarse gradient} of~$\calJ_n$
is given by $\nabla\calJ_n(m_n)=\calP_n^T\nabla\calJ(p_n)$, and hence, it is
simply obtained by summing the components of the fine gradient inside each zone
of the current parameterization~$\calP_n$.

A refinement of the current parameterization~$\calP_n$ is obtained by
{\em cutting} one zone~$Z_{n,j}$ into two subzones~$Z_{n,j_1}^\cut$
and~$Z_{n,j_2}^\cut$ with $1\leq j_1,j_2\leq n+1$.
The resulting $(n+1)$-zone zonation~$Z_n^\cut$ is a candidate for~$Z_{n+1}$.
This operation is fully characterized by the subsets of indices corresponding
to the new subzones which satisfy the disjoint union
$I_{n,j}=I_{n,j_1}^\cut\sqcup I_{n,j_2}^\cut$.
The associated parameterization is denoted by~$\calP_n^\cut$ and the
corresponding objective function by~$\calJ_n^\cut$.
Opening the new (vectorial) degree of freedom while keeping the same value in
the two subzones does not change the optimum of the objective function, and the
components of the gradient $\nabla\calJ_n^\cut((m_n^\cut)^{0,\opt})$ are all
equal to zero except for components~$j_1$ and~$j_2$ that are opposite vectors
(of dimension~$n_p$).
From the first equation in~(\ref{eqn:lagrcond}) and the form of matrix~$A$, we
deduce that the Lagrange multiplier~$(\lambda_n^\cut)^{0,\opt}$ is equal to the
component~$j_2$ of the previous gradient.
Hence, the refinement indicator~$\calI_n^\cut$ associated with the current
refinement can be computed easily, at almost no additional cost, from the
partial derivatives of the objective function with respect to the fine
parameter in the subzones by
\begin{equation}
  \label{eqn:indic}
  \calI_n^\cut
  = \left\| \sum_{i/K_i\subset Z_{n,j_1}^\cut}
  \frac{\partial\calJ}{\partial p_i}(p_n^\opt)
  \right\|
  = \left\| \sum_{i/K_i\subset Z_{n,j_2}^\cut}
  \frac{\partial\calJ}{\partial p_i}(p_n^\opt)
  \right\|.
\end{equation}
Thus, we can define a large number of tentative cuttings producing possible
refinements in each zone, and compute their corresponding refinement
indicators.

Since the refinement indicators give us only a first order information on the
decrease of the optimal objective function, one can select a set of cuttings
associated with some of the highest values of the indicators, and not only with
the largest one.
Then, the objective function is minimized for all the parameterizations
obtained by implementing each of these selected cuttings.
And finally, the cutting defining the next zonation is the one leading to the
largest decrease of the objective function.
In the general case, this minimization phase is very expensive.
It is thus important to keep the number of selected cuttings very small.

The initial value of the coarse parameter for these minimizations is obtained
from the previous optimal value by duplicating the value in the split zone,
i.e.
\begin{equation}
  \label{eqn:mnjinit}
  (m_{n,j_1}^\cut)^\init = (m_{n,j_2}^\cut)^\init = m_{n,j}^\opt
  \qquad\mbox{and}\qquad
  (m_{n,j^\prime}^\cut)^\init = m_{n,j^\prime}^\opt
  \mbox{ for }j^\prime\not=j_1,j_2.
\end{equation}
The parameterization given by~(\ref{eqn:P_h}) is injective, and we can define
its least-squares pseudoinverse by
$\tilde{\calP}_n=(\calP_n^T\calP_n)^{-1}\calP_n^T$.
It computes the mean value on each zone.
When considering the scalar product weighted by the measure of each cell of the
fine mesh, the previous pseudoinverse corresponds to the {\em projection} map
given by
\begin{equation}
  \begin{array}{rrcl}
    \tilde{\calP}_n: & (\R^{n_p})^{n_I} & \longrightarrow & (\R^{n_p})^n \\
    & p_n = (p_{n,i})_{i\in I} & \longmapsto & m_n = (m_{n,j})_{1\leq j\leq n}
    \displaystyle \mbox{ with } m_{n,j} =
    \frac{\sum_{i\in I_{n,j}} \meas(K_i) p_{n,i}}{\meas(Z_{n,j})}.
  \end{array}
\end{equation}
Then, (\ref{eqn:mnjinit})~rewrites in the more compact form
$(m_n^\cut)^\init=\tilde{\calP}_n^\cut p_n^\opt$.

At the end, the next $(n+1)$-zone parameterization map~$\calP_{n+1}$ has only
one (vectorial) degree of freedom more, but it produces a significant benefit
effect on the objective function.

\begin{remark}
  \label{rem:finemesh}
  The adaptive parameterization technique is independent of the fine
  discretization~$\calT_h$ of the domain~$\Omega$.
  In particular, it is valid for regular or irregular meshes, structured or
  unstructured meshes.
\end{remark}

\subsection{A generic algorithm for adaptive parameterization}
\label{sec:generic}

The following algorithm can be applied to the estimation of distributed
multidimensional parameters in any partial differential equation.
\begin{description}
\item[Initialization] \mbox{}
  \begin{enumerate}
    \setcounter{enumi}{-2}
  \item \label{item:init}
    Choose an initial $n_0$-zone parameterization~$\calP_{n_0}$ and an initial
    coarse parameter~$m_{n_0}^\init$.
  \item \label{item:initmin}
    Compute $m_{n_0}^\opt=\arg\min\calJ_{n_0}(m_{n_0})$ from~$m_{n_0}^\init$
    and $p_{n_0}^\opt=\calP_{n_0}m_{n_0}^\opt$.
  \end{enumerate}

\item[Iterations]
  For $n\geq n_0$, do until convergence:
  \begin{enumerate}
  \item \label{item:finegradient}
    Compute the fine gradient $g_n^\opt=\nabla_p\calJ(p_n^\opt)$.
  \item \label{item:cuttings}
    Choose a set of cuttings~$C_n$ defining parameterization candidates
    $(\calP_n^k)_{k\in C_n}$.
  \item \label{item:indicators}
    For all $k\in C_n$, compute the corresponding indicator~$\calI_n^k$
    from~$g_n^\opt$ using (\ref{eqn:indic}).
  \item \label{item:largest}
    Compute the largest indicator~$\calI_n^{\max}$.
  \item \label{item:select}
    Select a subset $C_n(\calI_n^{\max})\subset C_n$ of cuttings associated
    with the highest indicators.
  \item \label{item:minimize}
    For all $k\in C_n(\calI_n^{\max})$, \\
    compute
    $(m_n^k)^\opt=\arg\min\calJ_n^k(m_n^k)$
    from~$(m_n^k)^\init=\tilde{\calP}_n^kp_n^\opt$.
  \item \label{item:keepbest}
    Keep the best cutting
    $k_0=\arg\min_{k\in C_n(\calI_n^{\max})}\calJ_n^k((m_n^k)^\opt)$ \\
    and set $\calP_{n+1}=\calP_n^{k_0}$, $m_{n+1}^\opt=(m_n^{k_0})^\opt$
    and $p_{n+1}^\opt=\calP_{n+1}m_{n+1}^\opt$.
  \end{enumerate}
\end{description}

This generic algorithm may be adapted by the user to match his/her needs.
The convergence criterion can be specified, for instance: a small value of the
objective function (meaning that the data are correctly fitted), a small value
of the largest refinement indicator (meaning that there is no more important
discontinuities left to discover), the maximum number of zones is reached\ldots
In addition, the user is free to choose the refinement strategy defining the
full set of cuttings and the subset of those with the highest indicators.

Furthermore, as in~\cite{benchajaf}, it is possible to add at the end of each
iteration a coarsening step that allows zones corresponding to coarse
parameters with close values to be merged.

\subsection{Best cutting for multidimensional parameters}
\label{sec:best}

Let us suppose that we want to find among all possible cuttings the one
corresponding to the best indicator.
For example, this is useful in section~\ref{sec:segmentation} with the
identity direct operator for which this best cutting is related to the best
decrease of the objective function.
But this {\em absolute} best cutting may also be useful in the general case to
quantify the adequacy of the chosen set of cuttings to the inverse problem
trough the ratio of the best indicator within the chosen set of cuttings to the
absolute best indicator.
A high ratio shows that there exists cuttings outside the chosen set that
provide higher indicators.

In the scalar case, it is obvious from equation~(\ref{eqn:indic}) that the best
cutting for a given zone corresponds to follow the sign of the fine
gradient~(\ref{eqn:GpJn}).
In the multidimensional case, this is no more so simple as finding the best
cutting amounts to solve a discrete optimization problem.
We give here several possible heuristics.

The simpler idea consists in defining the sign of a vector and then to use the
same technique as in the scalar case.
There are several possible definitions for the sign of a vector.
Of course, none of them are fully satisfactory, but they can provide pertinent
zonations, as it can be seen in section~\ref{sec:segmentation} where the sign
of a three-dimensional vector is obtained by {\em majority vote} of the signs of
all components.
It could also have been defined as the sign of the sum of all components of the
vector.

Another simple idea is to consider a zonation for each component of the
parameter.
Then, the algorithm for an $n_p$-dimensional parameter on a grid of size~$n_I$
would be equivalent to the algorithm for a scalar parameter on a grid of size
$n_p\times n_I$.
Again, with the example of section~\ref{sec:segmentation} where the data are
colored images, $n_p=3$ and~$n_I$ is the number of pixels.
In this case, the algorithm for an RGB (Red/\-Green/\-Blue) image of size~$n_I$
would be equivalent to the algorithm for 3~grayscale images of the same
size~$n_I$.
This kind of decomposition reminds the dawn of color photography where three
black-and-white pictures were taken with three colored filters, and more
recently the 3-LCD video projectors.
The main drawback of this technique is that the compound zonation corresponding
to the superimposition of the zonations for all the components may contain a
lot of small zones.
Indeed, the intersection of~$n_p$ zonations with~$n$ zones each can produce up
to~$n^{n_p}$ zones.
In other words, we would add up to $\sum_{k=0}^{n_p-1}C_{n_p}^kn^k$ degrees of
freedom per component at each iteration!
Which is completely in contradiction with the adaptive parameterization goal of
providing the coarsest parameter, i.e. adding at most only one degree of
freedom per component per iteration.

Keeping the attractive idea of computing a {\em local} indicator for
each---scalar---component of the parameter separately%
\footnote{For which the {\em best} cutting associated with the highest local
  indicator simply follows the sign of the---scalar---fine gradient for each
  component.}, we can select at each iteration among the best cuttings for each
component the one providing the highest {\em global} indicator.
Then, this cutting is applied to all components, thus guaranteeing the addition
of a single degree of freedom per component.
Of course, this heuristics does not always provide the highest possible global
indicator, but it seems reasonable to us and we suggest to use it.

\subsection{Implementation issues}
\label{sec:implementation}

The first issue is the genericity of the implementation of the algorithm.
Numerical results shown in~\cite{benchajaf} were obtained with a Fortran~77
program where the refinement algorithm was hard-coded within the model for the
very specific problem of the estimation of hydraulic transmissivities in a
parabolic flow equation set on a planar regular mesh.
Once the technique has been validated on this example, it is very advisable to
propose to the community a generic implementation that will allow to apply the
refinement indicator adaptive parameterization algorithm to any problem
implemented in a wide variety of programming languages.

The second issue is related to the definition of the set of cuttings, and more
particularly to the implementation of it.
In the general case, it is not possible to compute the true parameterization,
and then, the adaptive parameterization technique is based on the definition of
a set of cuttings, which would preferably be adapted to the problem. 
In~\cite{benchajaf}, the two-dimensional mesh was supposed regular and
rectangular, and the refinements were of four elementary kinds: horizontal
cutting, vertical cutting, oblique cutting and checkerboard cutting.
At each iteration, the indicators related to all possible cuttings within the
selected family for all current zones were computed.
Testing larger sets of cuttings, e.g. by enriching the elementary cutting
offer, is very interesting because it will allow the user to use a priori
information adapted to the problem.
Furthermore, the ability to specify that some elementary cuttings should only
apply to a particular part of the domain will allow to interpret more closely 
the a priori information.

The last issue concerns the performance of the implemented code.
At each iteration, steps~\ref{item:indicators} and~\ref{item:minimize} of the
algorithm correspond respectively to the computation of all the indicators
associated with the chosen cuttings and to the actual minimization for the
selected cuttings.
Both steps involve fully independent computations that may run in parallel.

All these aspects show the necessity for high-level programming capabilities.
In~\cite{benclewei3}, we choose Objective Caml (\ocaml) which is a fast
modern high-level functional programming language based on sound theoretical
foundation and formal semantics%
\footnote{See~{\tt http://www.ocaml.org/}.}.
It is particularly well suited for the implementation of complex algorithms.
The generic driver is written in \ocaml\/ and data exchange with the worker
follow a simple and safe communication protocol that is implemented in several
common programming languages including \ocaml, \clang, \cpp\/ and \fortran.
Here the generic tasks sent to the worker are {\tt "cost"}, {\tt "grad"} and
{\tt "optim"} respectively to implement the functions $p\longmapsto\calF(p)$,
$p\longmapsto\nabla\calJ(p)$ and
$(m,\calP)\longmapsto(m^\opt,\calJ^\opt)
=(\arg\min\calJ_\calP(m),\calJ_\calP(m^\opt))$.
The design of a mini-language allow a flexible and convenient way to define the
set of cuttings.
Furthermore, {\em automatic} parallelization capabilities through the skeleton
programming system provided by \ocamlpppl%
\footnote{See~{\tt http://www.ocamlp3l.org/.}}
has recently been successfully used in the field of Scientific Computation for
domain decomposition applications, see~\cite{fcl:cp:05,fcl:parco:06}.

% $Id: segmentation.tex,v 1.2 2008/01/16 12:32:01 fclement Exp $

\section{Application to the identity model}
\label{sec:segmentation}

To illustrate the use of our \ocaml\/ implementation, detailed
in~\cite{benclewei3}, and the capabilities of the refinement indicators
algorithm in the multidimensional case, presented in the previous section, we
consider the simplest model ever: the identity model.
Thus, in this section, we assume that $\calF=\Id$.
Of course, in this very simple case, the inverse problem does not present any
difficulty as the minimum of the objective function~(\ref{eqn:J_h}) is exactly
zero, and it is obtained for $p^\opt=d$.
Note that the fine gradient is then simply given by the difference
$\nabla\calJ(p)=p-d$.
Nevertheless, assuming that the data is actually an image of size
$n_I=n_x\times n_y$, then the mesh is regular and rectangular, the cells are
called {\em pixels}, and the application of the refinement indicators algorithm
amounts in this case to determine groups of pixels with the same color---or
gray level.
This is the basic idea of image segmentation for edge detection, shape
recognition, compression/simplification, and so on\ldots

The segmentation of a grayscale image is simply a scalar example, but even in
this case, the determination of the best cutting providing the highest decrease
of the objective function is not an easy task, see~\cite{benclewei2}.
The case of color images is more complicated as the color information is
typically multidimensional.
In computer applications, it is usually represented using a three-dimensional
color space---hence $n_p=3$: either RGB, for the Red/\-Green/\-Blue color cube,
HSV, for the Hue/\-Saturation/\-Value color cone, or HSL, for the
Hue/\-Saturation/\-Lightness color double cone.
Of course, it is possible to define a scalar discrete representation of colors,
e.g. by packing the 24-bit RGB encoding into a single integer ranging from~0
to~$2^{24}-1$.
But, such a nonlinear mapping produces weird effects when averaging on a zone:
the mean color between red pixels and blue pixels may become an odd green,
instead of the natural purple value.
We consider here the RGB representation of color images.

In the case of the identity model, the operator $\calF_n^T\calF_n$ of
subsection~\ref{sec:quantif} reduces to the diagonal matrix $\calP_n^T\calP_n$
which gives the numbers of pixels per subzone.
Then, equation~(\ref{eqn:decrease}) rewrites
\begin{equation}
  \label{eqn:decrease_id}
  \calJ_1^\opt - \calJ_2^\opt =
  \frac{1}{2} \frac{n_j}{n_{j_1}n_{j_2}} \calI_{n,j}^2
\end{equation}
where~$n_j$ is the number of cells in the selected zone~$Z_{n,j}$, and $n_{j_1}$
and~$n_{j_2}$ are respectively the number of cells in the
subzones~$Z_{n,j_1}^\cut$ and~$Z_{n,j_2}^\cut$.
Although the decrease of the objective function provided by a cutting is
proportional to the square of the refinement indicator associated with this
cutting, unfortunately the proportionality coefficient depends on the cutting
itself.
Nevertheless, for high values of~$a$, the function
$x\longmapsto\frac{a}{x(a-x)}$ is almost flat in the neighborhood of~$a/2$, and
the highest indicator almost yields the best decrease of the objective
function, at least during the first iterations of the algorithm.
Furthermore, there is no more interaction between the pixels during the
minimization process.
Hence, there is no need to compute the refinement indicators for all zones at
each iteration, but rather, it is advisable to compute only the indicators
related to the two most recent subzones, and to keep in memory the values for
the unchanged zones.
An optimized version of the refinement indicators algorithm dedicated to the
diagonal models is described in~\cite{benclewei2}.

For this first multidimensional application of the algorithm, we have simply
defined the {\em pseudo}-sign of the fine vector gradient on each pixel by the
sum of the signs for each component.
Let $\sgn$ be the sign function on the field of real numbers (i.e.,
$\sgn(x)=-1$ when $x<0$, $\sgn(x)=+1$ when $x>0$ and $\sgn(0)=0$).
Then, we define the sign of any triple by
$\sgn(x,y,z)=\sgn(\sgn(x)+\sgn(y)+\sgn(z))$.
As developed in subsection~\ref{sec:best}, this is a very basic heuristics to
deal with the multidimensional case, but it is sufficient for our present
illustration purposes.

In subsection~\ref{sec:generic}, the choice for the refinement strategy
defining the set of cuttings and the subset of those with the highest
indicators was let to the user.
Here, we consider two strategies.
The first one is a mild flavor of the {\em best} strategy consisting in
choosing at each iteration the cutting related to the highest indicator which
guarantees in the case of the identity model the best descent of the objective
function.
In fact, we follow the pseudo-sign of the gradient as defined above.
This strategy actually provides an image segmentation technique.
The second one is called the {\em dichotomy} strategy.
It is an example of strategy based on elementary cuttings deployed in all
zones.
The only cuttings allowed are those splitting a zone vertically or horizontally
in its middle.
This illustrates what could happen in the general case where the best strategy
is not necessarily optimal and where the adaptive parameterization technique is
based on the choice for a set of cuttings.
Notice that this could provide interesting creative filters for image
processing.

The next Figures follow the same scheme.
They are all made of four images.
The data image~$d$ is displayed on the lower left corner.
The segmented image~$p_n^\opt$ is displayed---in true colors---on the lower
right corner, and the corresponding $n$-zone zonation~$\calP_n$ is displayed in
false colors on the upper right corner.
The pseudo-sign of the fine gradient $\nabla\calJ(p_n^\opt)$ is displayed on
the upper left corner where positive values are depicted in white, negative
values in black and small absolute values in gray.
Some monitoring counters are also displayed.
This is actually the interactive display of the \ocaml\/ driver at runtime,
see~\cite{benclewei3}.

\bigskip

\begin{figure}[p]
  \begin{center}
    \includegraphics[width=10cm]{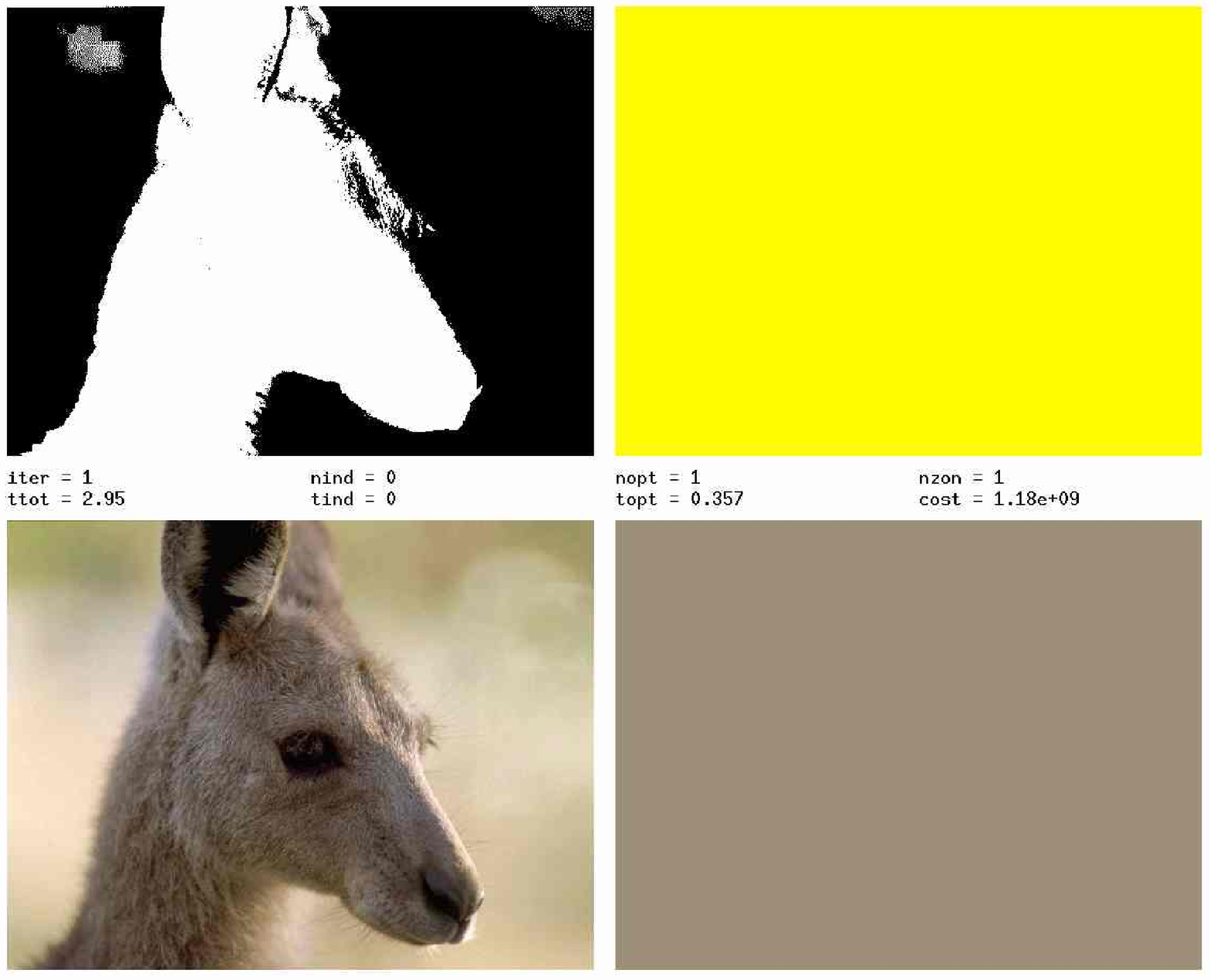}
    \caption{Best strategy, 1 zone.
      Bottom: data image~$d$ (left), segmented image~$p_1^\opt$ (right).
      Top: pseudo-sign of the gradient~$g_1^\opt$ (left), 1-zone
      zonation~$\calP_1$ (right).}
    \label{fig:best01}
  \end{center}
\end{figure}

\begin{figure}[ht]
  \begin{center}
    \includegraphics[width=10cm]{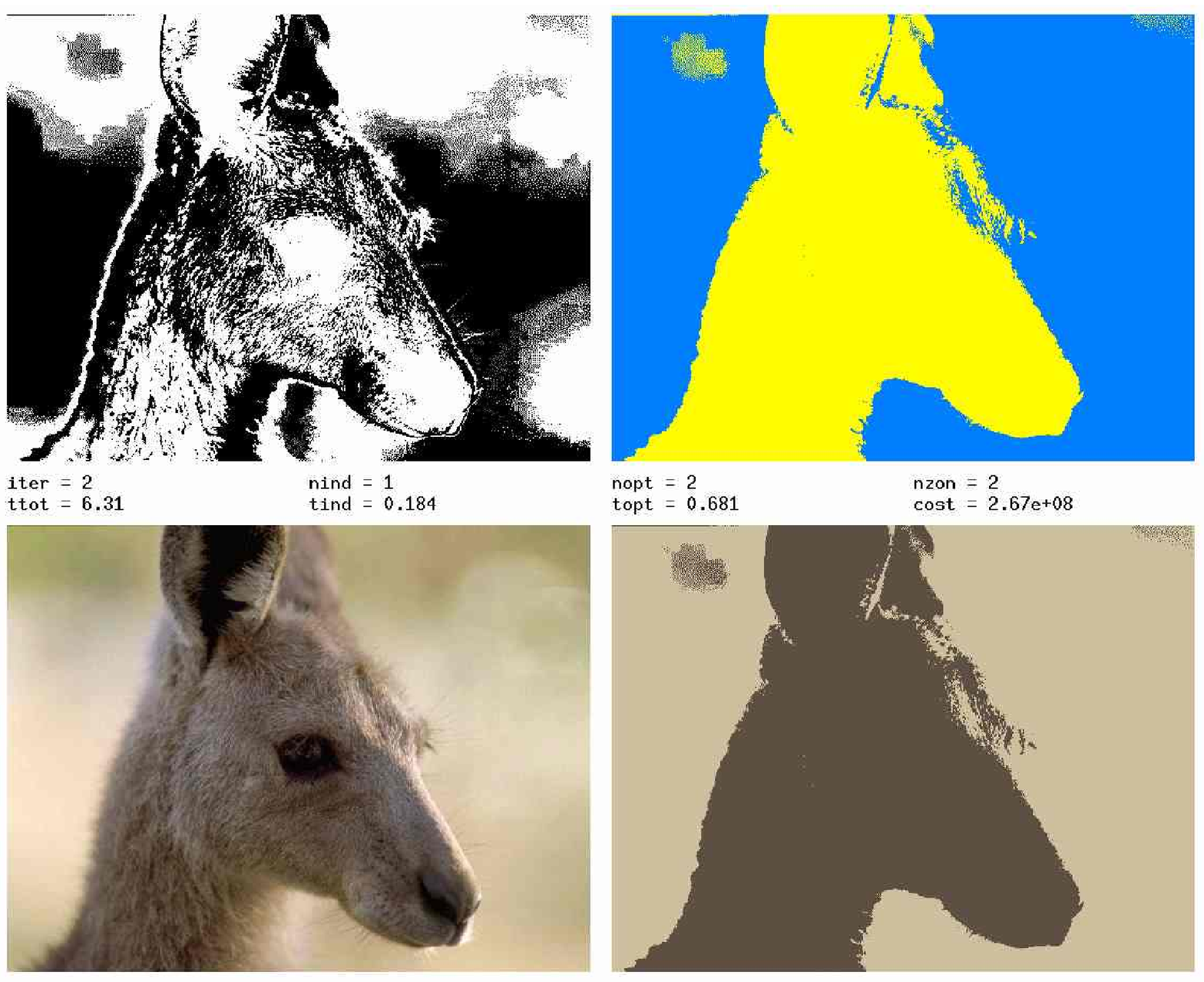}
    \caption{Best strategy, 2 zones.
      Bottom: data image~$d$ (left), segmented image~$p_2^\opt$ (right).
      Top: pseudo-sign of the gradient~$g_2^\opt$ (left), 2-zone
      zonation~$\calP_2$ (right).
      The refinement indicator for the first cutting was
      $\calI_1^\best=1.0\,10^7$.}
    \label{fig:best02}
  \end{center}
\end{figure}

\begin{figure}[ht]
  \begin{center}
    \includegraphics[width=10cm]{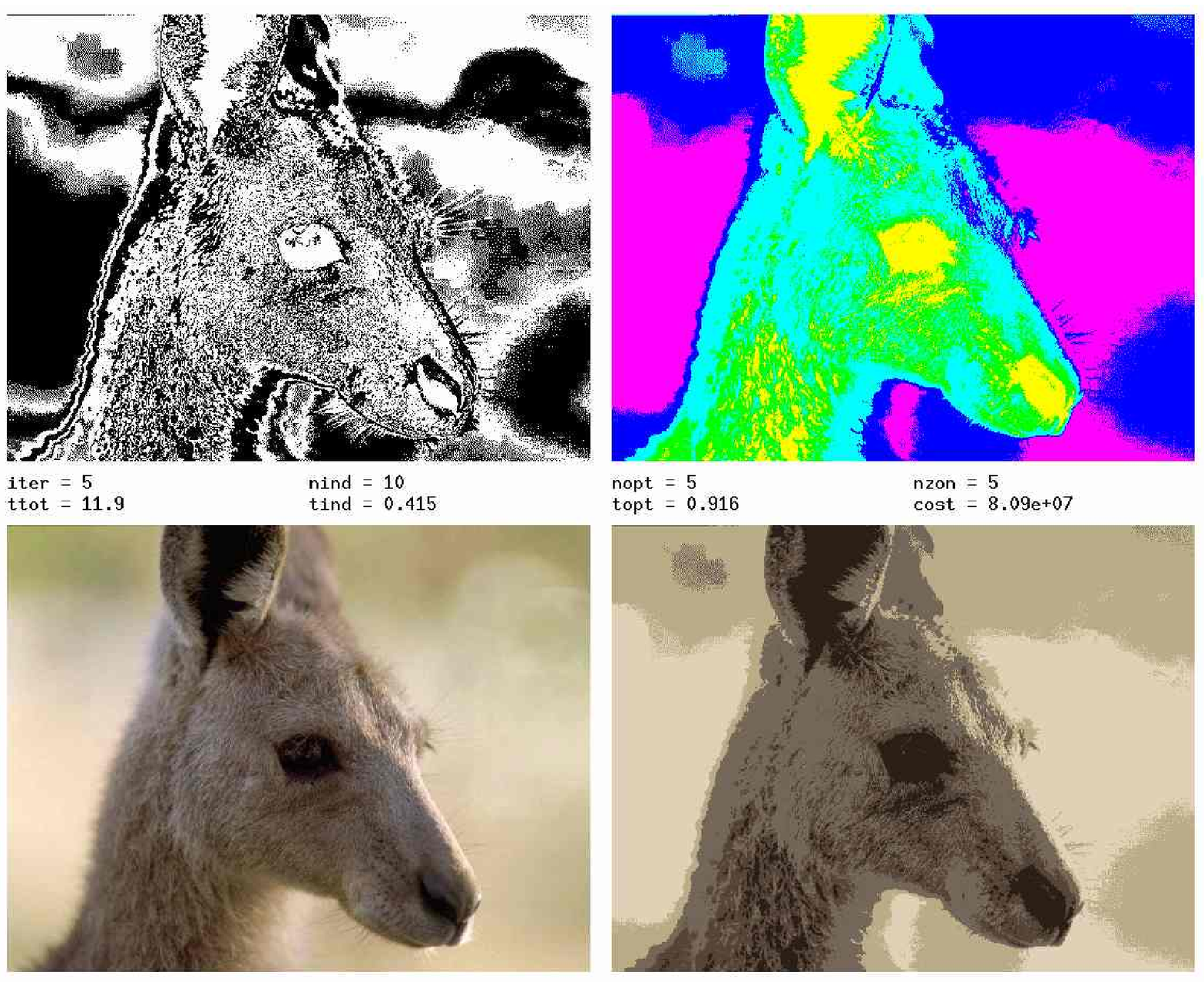}
    \caption{Best strategy, 5 zones.
      Bottom: data image~$d$ (left), segmented image~$p_5^\opt$ (right).
      Top: pseudo-sign of the gradient~$g_5^\opt$ (left), 5-zone
      zonation~$\calP_5$ (right).
      The refinement indicator for the fourth cutting was
      $\calI_4^\best=7.5\,10^5$.}
    \label{fig:best05}
  \end{center}
\end{figure}

\begin{figure}[ht]
  \begin{center}
    \includegraphics[width=10cm]{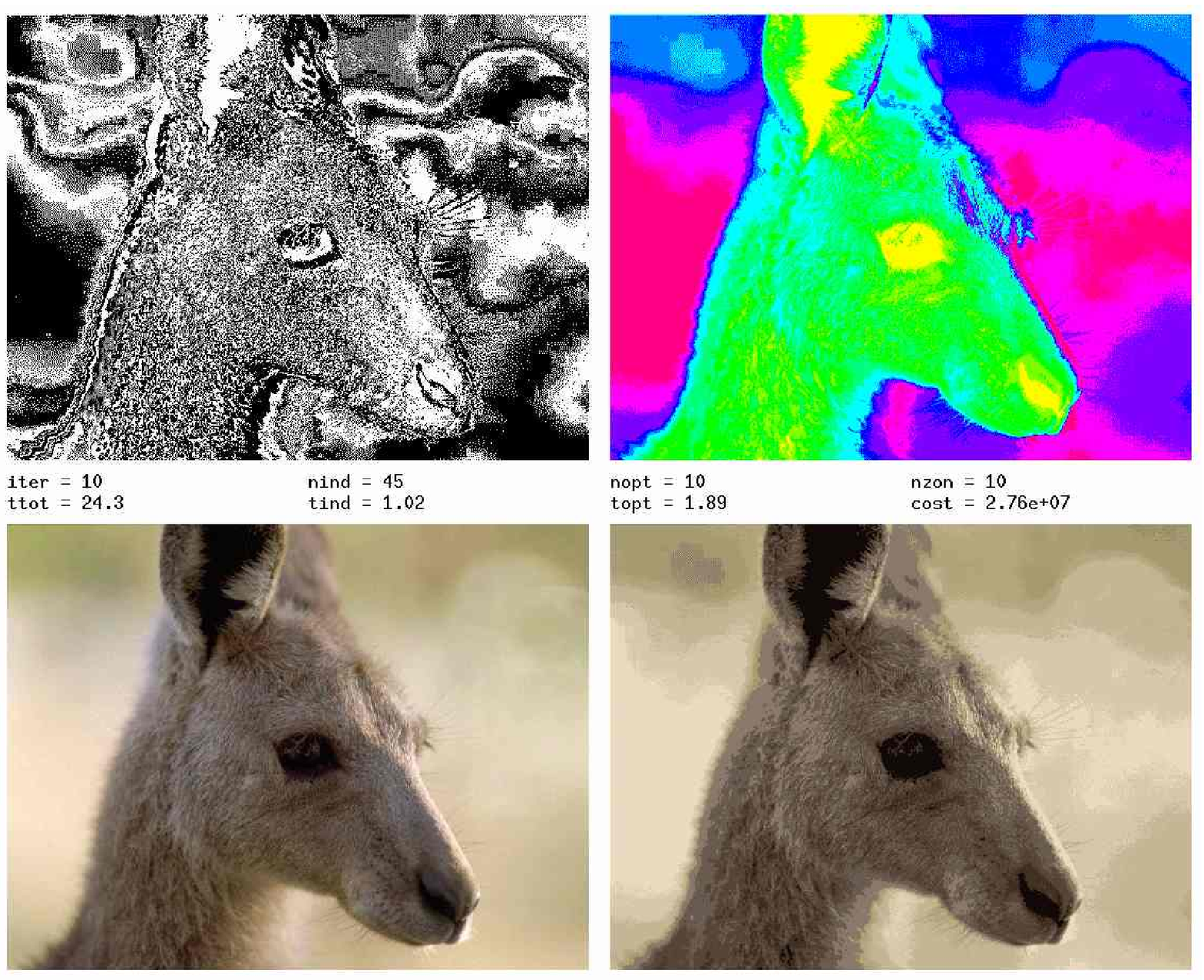}
    \caption{Best strategy, 10 zones.
      Bottom: data image~$d$ (left), segmented image~$p_{10}^\opt$ (right).
      Top: pseudo-sign of the gradient~$g_{10}^\opt$ (left), 10-zone
      zonation~$\calP_{10}$ (right).
      The refinement indicator for the ninth cutting was
      $\calI_9^\best=1.7\,10^5$.}
    \label{fig:best10}
  \end{center}
\end{figure}

%iter = 2: zone = 1 -> indic = 1.01327e+07, 
%iter = 3: zone = 2 -> indic = 2.28113e+06, 
%iter = 4: zone = 1 -> indic = 1.99997e+06, 
%iter = 5: zone = 1 -> indic = 754217, 
%iter = 6: zone = 4 -> indic = 728514, 
%iter = 7: zone = 3 -> indic = 619458, 
%iter = 8: zone = 7 -> indic = 535370, 
%iter = 9: zone = 1 -> indic = 297445, 
%iter = 10: zone = 6 -> indic = 167972, 

The first series of experiments exploit the best strategy, hence we may expect
the fastest decrease of the objective function.
The results obtained after iterations number~1, 2, 5 and~10 are respectively
represented in Figure~\ref{fig:best01}, \ref{fig:best02}, \ref{fig:best05}
and~\ref{fig:best10}.

At the beginning, in Figure~\ref{fig:best01}, the initial zonation contains only
one zone, the initial parameter is $m_1^\init=(0,0,0)$---i.e., a black
image---and the initial value of the objective function is
$\calJ_1^\init=8.0\,10^9$.
After the initialization phase, the objective function is
$\calJ_1^\opt=1.2\,10^9$, meaning that 62\% of the data are explained%
\footnote{At iteration~$n$, the percentage of data explained is defined as
  $100\times\left(1-\sqrt{\frac{\calJ_n^\opt}{\calJ_1^\init}}\right)$.}.
Then, the pseudo-sign of the first fine gradient $g_1^\opt$ already depicts the
head of the kangaroo.
This cutting is associated with the---high---refinement indicator
$\calI_1^\best=1.0\,10^7$.

After the fifth iteration, in Figure~\ref{fig:best05}, 90\% of the data are
explained and the main features of the original picture are visible.
The last highest indicator was already only 7.4\% of the first one.
With 10~zones, in Figure~\ref{fig:best10}, 94\% of the data are explained.
For the eye, the segmented image is very close to the original one.
The last highest indicator was a flat 1.7\% of the first one.

\bigskip

\begin{figure}[ht]
  \begin{center}
    \includegraphics[width=10cm]{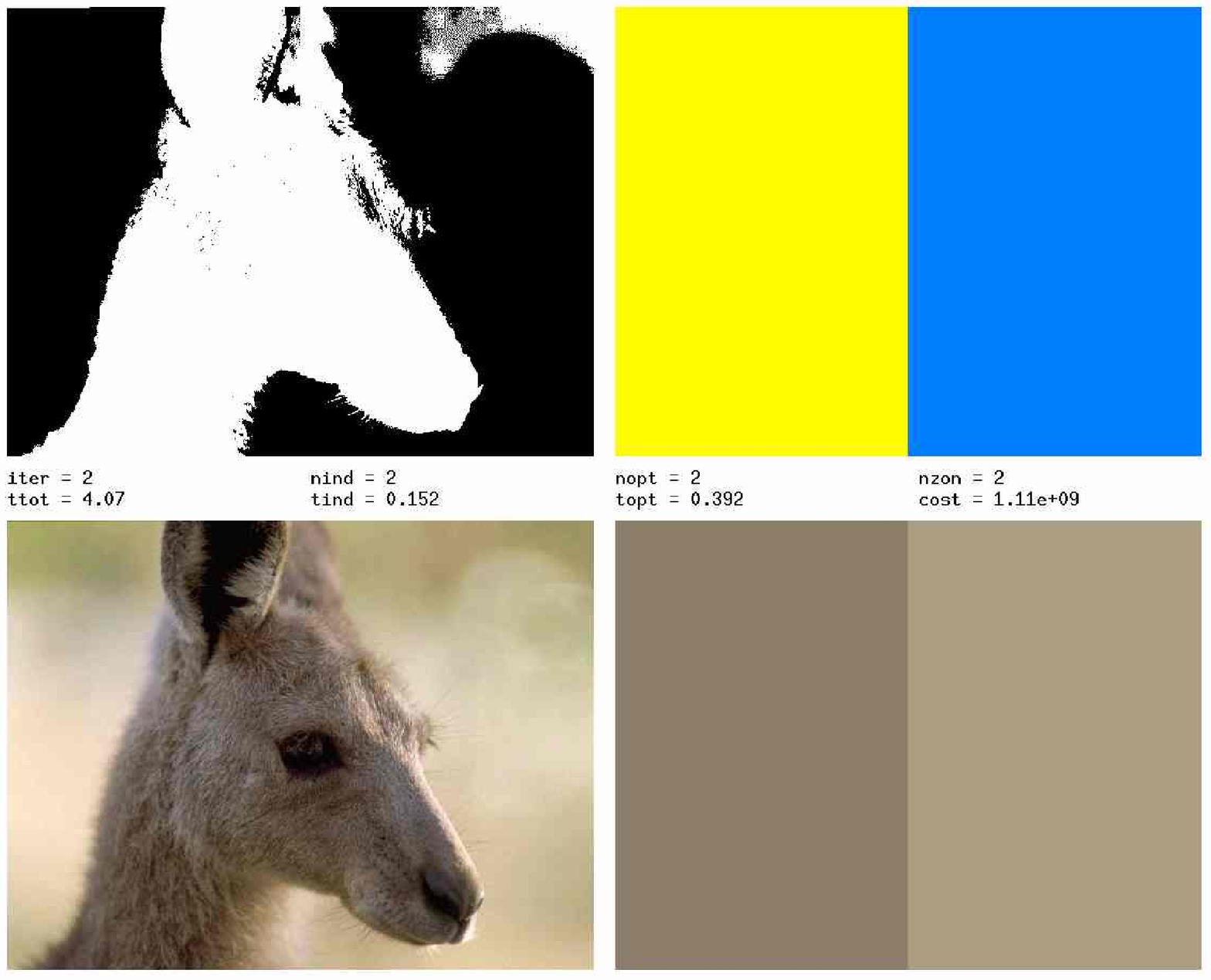}
    \caption{Dichotomy strategy, 2 zones.
      Bottom: data image~$d$ (left), segmented image~$p_2^\opt$ (right).
      Top: pseudo-sign of the gradient~$g_2^\opt$ (left), 2-zone
      zonation~$\calP_2$ (right).
      The refinement indicator for the first cutting was
      $\calI_1^\dichotomy=2.9\,10^6$.}
    \label{fig:dichotomy02}
  \end{center}
\end{figure}

\begin{figure}[ht]
  \begin{center}
    \includegraphics[width=10cm]{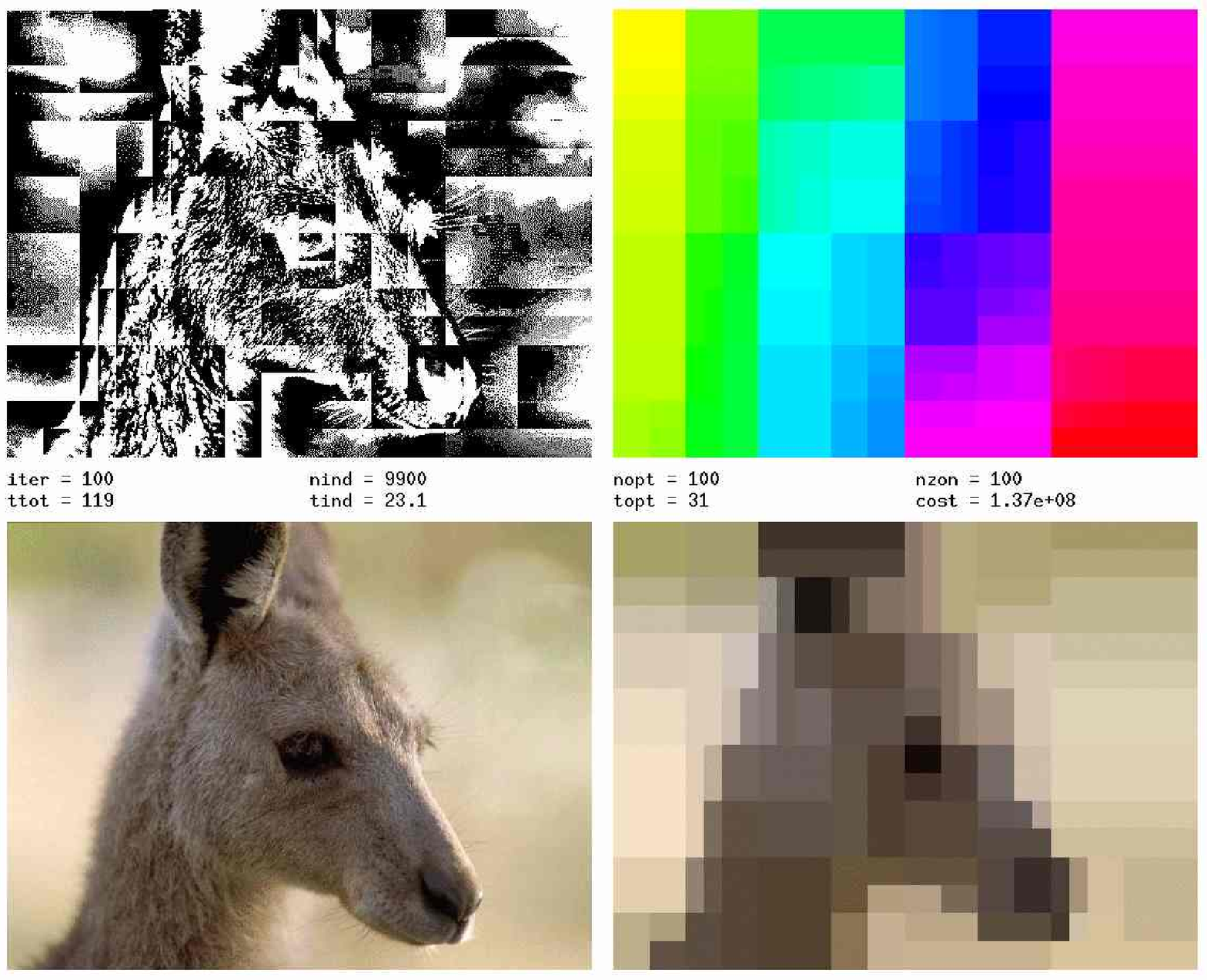}
    \caption{Dichotomy strategy, 100 zones.
      Bottom: data image~$d$ (left), segmented image~$p_{100}^\opt$ (right).
      Top: pseudo-sign of the gradient~$g_{100}^\opt$ (left), 100-zone
      zonation~$\calP_{100}$ (right).
      The refinement indicator for the 99th cutting was
      $\calI_{99}^\dichotomy=2.8\,10^4$.}
    \label{fig:dichotomy100}
  \end{center}
\end{figure}

%iter = 2: zone = 1 -> indic = 2.85802e+06, 
%iter = 3: zone = 1 -> indic = 4.0436e+06, 
%iter = 4: zone = 3 -> indic = 2.80378e+06, 
%iter = 5: zone = 1 -> indic = 1.27809e+06, 
%iter = 6: zone = 4 -> indic = 1.11259e+06, 
%iter = 7: zone = 2 -> indic = 818940, 
%iter = 8: zone = 5 -> indic = 553533, 
%iter = 9: zone = 3 -> indic = 321042, 
%iter = 10: zone = 8 -> indic = 312501, 
%iter = 30: zone = 26 -> indic = 144632, 
%iter = 100: zone = 72 -> indic = 27657.1, 

The second series of experiments use the dichotomy strategy.
The result obtained after~2 and~100 iterations are respectively represented in
Figures~\ref{fig:dichotomy02} and~\ref{fig:dichotomy100}.

The first iteration is exactly the same as with the best strategy in
Figure~\ref{fig:best01}.
Then, the first cutting does not reveal at all the head of the kangaroo.
The image with two zones, in Figure~\ref{fig:dichotomy02}, only explain one
percent more of the data than the one with one zone.
The corresponding indicator is only 28\% of the first indicator with the best
strategy.
Indeed, the chosen family of cuttings is not well suited for the problem.
Even at iteration number~100, in Figure~\ref{fig:dichotomy100}, only 87\% of
the data are explained and the last highest indicator was about 1\% of the
first one.

Of course, this strategy is inefficient for the segmentation of images, but
exploring a subset of possible cuttings at each iteration acts as a
preconditioner, and it may also help in avoiding to fall into local minima in
the general case.
Notice that after 100~iterations of the dichotomy strategy, the gradient keeps
very sharp details that has already disappeared after only 10`iterations of the
best strategy.
We could have done a better job by allowing any vertical and horizontal
cuttings, but there would have been much more indicators to evaluate.
The dichotomy run took 54~seconds, compared to the 3~seconds for the best run
(only add {\tt tind} and {\tt topt}, as the total time {\tt ttot} also takes
into account the display).

% $Id: conclusions.tex,v 1.2 2008/01/16 12:32:01 fclement Exp $

\section{Conclusions}

We give a general formulation of the refinement indicator algorithm for the
estimation of \emph{vector valued} distributed parameters in \emph{any} partial
differential equation. \\
The main findings are: \\
(i) in the multidimensional case, the refinement indicator associated with a new
degree of freedom is the norm of the corresponding Lagrange multiplier. \\
(ii) In the linear case, the refinement indicator also measures the decrease
of the optimum least-squares data misfit objective function when the
corresponding degree of freedom is open. \\
(iii) When applied to the inversion of the identity, the vector valued version
of the refinement indicators algorithm provides a surprisingly powerful tool
for the segmentation of color images (where colors are three-dimensional
vectors in the RGB model).

% $Id: acknowledgments.tex,v 1.2 2008/01/16 12:32:01 fclement Exp $

\section*{Acknowledgments}

The authors thank the GDR MOMAS for providing partial financial support.

\bibliography{biblio}
\bibliographystyle{plain}

\end{document}